\documentclass[preprint,number]{elsarticle_mod}
\usepackage{graphicx}

\usepackage{amsmath}
\usepackage{amsfonts}
\usepackage{amssymb}
\usepackage{amsbsy} 

\usepackage{amsthm}
\newtheorem{thm}{Theorem}
\newtheorem{conj}{Conjecture}

\theoremstyle{definition} 
\theoremstyle{theorem} 
\theoremstyle{theorem} \newtheorem*{remark}{Remark}
\theoremstyle{definition} \newtheorem{example}{Example}

\def\bs#1{\boldsymbol#1}

\newcommand{\edvu}{v \rightarrow u}
\newcommand{\eduv}{u \rightarrow v}

\begin{document}
\begin{frontmatter}



\title{A conjecture on independent sets and graph covers}


\author[ism]{Yusuke Watanabe}
\ead[ism]{watay@ism.ac.jp}



\address[ism]{The Institute of Statistical Mathematics 
10-3 Midori-cho, Tachikawa, Tokyo 190-8562, Japan
Phone: +81-(0)50-5533-8500}

\begin{abstract}
In this article, I present a simple conjecture on the number of independent sets on graph covers.
The conjecture implies that the partition function of a binary pairwise attractive model
is greater than that of the Bethe approximation.
\end{abstract}

\begin{keyword}
graph cover \sep
independent set \sep
Bethe approximation \sep
\end{keyword}

\end{frontmatter}

\section{Terminologies}
Throughout this article, $G=(V,E)$ is a finite graph with vertices $V$ and undirected edges $E$.
For each undirected edge of G, we make a pair of oppositely directed edges, which form a set of
directed edges $\vec{E}$. Thus, $|\vec{E}|= 2 |E|$.
%

An {\it $M$-cover} of a graph $G$ is its $M$-fold covering space\footnote{Interpret graphs as topological spaces.}. 
All $M$-covers are explicitly constructed using permutation voltage assignment as follows \cite{gross1977generating}.
A {\it permutation voltage assignment} of $G$ is a map
\begin{equation}
 \alpha : \vec{E} \rightarrow \mathfrak{S}_M  \qquad \text{s.t.} \quad \alpha(\eduv)=\alpha(\edvu)^{-1} \quad \forall uv \in E,
\end{equation}
where $\mathfrak{S}_M$ is the permutation group of $\{1,\dots,M\}$.
Then an $M$-cover $\tilde{G}=(\tilde{V},\tilde{E})$ of $G$ is given by 
$\tilde{V}:=V \times \{1,\dots,M\}$ and 
\begin{equation}
 (v,k)(u,l) \in \tilde{E} \Leftrightarrow uv \in E \text{ and } l = \alpha(\edvu)(k).
\end{equation}
If an $M$-cover is $M$ copies of $G$ then it is called {\it trivial $M$-cover} and denoted by $G^{\oplus M}$. 
This is obtained by identity permutations.
The {\it natural projection}, $\pi$, from a cover $\tilde{G}$ to $G$ is obtained by forgetting the ``layer number''.
That is, $\pi : \tilde{V} \rightarrow V$ is given by $\pi(u,i)=u$ and $\pi : \tilde{E} \rightarrow E$ is given by $\pi((v,k)(u,l))=vu$. 

An independent set $I$ of a graph $G$ is a subset of $V$ such that none of the elements in $I$ are adjacent in $G$.
Formally, $I$ is an independent set iff $u,v \in I \Rightarrow uv \not\in E$.
The multivariate {\it independent set polynomial} of $G$ is defined by
\begin{equation}
 p(G):= \sum_{I \text{: independent set}}  \quad  \prod_{v \in I} x_v,
\end{equation}
with indeterminates $x_v \quad (v \in V)$.

\section{The conjecture}
We extend the definition of the projection map $\pi$ over the multivariate polynomial
ring.
First, let us define a map $\Pi$ for each indeterminate by
$\Pi(x_v)=x_{\Pi(v)}$, where $v \in \tilde{V}$ . 
Then this is uniquely extended to the polynomial ring as a ring
homomorphism; for example $\Pi(x_v + x_{v'})= \Pi(x_v)+\Pi(x_{v'})$
and $\Pi(x_v  x_{v'})= \Pi(x_v)\Pi(x_{v'})$.
\begin{conj}
\footnote{I have checked the conjecture for many examples by computer.}
\label{conj:indep}
 For any bipartite graph $G$ and its $M$-cover $\tilde{G}$, we
 conjecture the following relation:
\begin{equation}
 \Pi( p(\tilde{G} ) )  \preceq p(G)^M,
\end{equation}
where $\Pi$ is defined as above and the symbol $\preceq$ means the
inequalities for all coefficients of monomials.
\end{conj}
We can interpret the conjecture more explicitly as follows.
For a subset $U$ of $V$, define 
$\mathcal{I}(\tilde{G},U):=\{I \subset \tilde{V} |I \text{ is
independent set}, \pi(I)=U \}$, 
where $\pi (I)$ is the image of $I$.
Since $p(G)^M = \Pi( p(G^{\oplus M}))$, the conjecture is equivalent to the following statement:
\begin{equation}
 |\mathcal{I}(\tilde{G},U)| \leq  |\mathcal{I}(G^{\oplus M},U)|  \text{
  for all }  U \subset V. \label{eq:Zfinite}
\end{equation}
\begin{example}
Let $G$ be a cycle graph of length four and let $\tilde{G}$ be its 3-cover that is isomorphic to the cycle of length twelve.
Then,
\begin{equation}
 p(G)=1+ x_1 + x_2 + x_3 + x_4 + x_1 x_3 + x_2 x_4,
\end{equation}
\begin{equation}
 p(\tilde{G}) = 1+ \sum_{v=1}^{4} \sum_{m=1}^3 x_{(v,m)} + \ldots ,
\end{equation}
\begin{equation}
 \Pi(  p(\tilde{G}) ) =  1 + 3(x_1 + x_2 + x_3 + x_4)+  \ldots.
\end{equation}
It takes time and effort to check the conjecture, however, it is true in this case.
\end{example}

\begin{remark}
{\rm
The above conjecture is claimed for the pair {\tt (bipartite graph, independent set)}.
I also conjecture analogous properties for 
{\tt (bipartite graph, matching)},
{\tt (graph with even number of vertices, perfect matching)}
and
{\tt (graph, Eulerian set\footnote{A subset of edges is Eulerian if it induces a subgraph that only has vertices of degree two and zero.})}.
}
\end{remark}

\section{Implication of the conjecture}
The conjecture originates from the theory of the Bethe approximation.
The {\it partition function} of a binary pairwise model on a graph $G$ is 
\begin{equation}
 Z(G;\bs{J},\bs{h}) := \sum_{\bs{s} \in \{0,1\}^V}
\exp (\sum_{uv \in E} J_{uv} s_u s_v + \sum_{v \in V} h_v s_v),
\end{equation}
where the weights $(\bs{J},\bs{h})$ are called {\it interactions}.\footnote{In the following,
for a cover $\tilde{G}$ of $G$, we think that interactions are naturally induced from $G$.} 
It is called {\it attractive} if $J_{vu} \geq 0$ for all $vu \in E$.
The {\it Bethe partition function}\footnote{This is computed from the absolute minimum of the Bethe free energy; other Bethe approximations of the partition function corresponding to local minima are smaller than $Z_B$.}
$Z_B$ is defined by \cite{vontobelcounting}
\begin{align}
 Z_B :&= \exp \Big( - \min_{q} F_B(q) \Big)  \\
 &=\limsup_{M \rightarrow \infty} {< Z(\tilde{G})>}^{1/M},  \label{eq:PascalChar}
\end{align}
where $F_B$ is the Bethe free energy and $< \cdot >$ is the
mean with respect to the $M!^{|E|}$ covers.
(Details are omitted. See \cite{vontobelcounting}.)

\begin{thm}
If Conjecture \ref{conj:indep} holds, then
\begin{equation}
  Z \geq Z_B
\end{equation}
holds for any binary pairwise attractive models.\footnote{In a quite limited situation, the inequality is proved in \cite{sudderth2008loop}.}
\end{thm}
\begin{proof}
 From (\ref{eq:PascalChar}), the assertion of the theorem is proved if we show that 
\begin{equation}
 Z(G)^M \geq Z(\tilde{G})
\end{equation}
for any $M$-cover $\tilde{G}$ of $G$.
In the following, we see that the partition function can be written by the independent set polynomial and thus
the above inequality holds under the assumption of Conjecture \ref{conj:indep}.

\begin{align}
 Z(G)&=\sum_{\bs{s} \in \{0,1\}^V}
\exp (\sum_{uv \in E} J_{uv} s_u s_v + \sum_{v \in V} h_v s_v) \\
& = \sum_{\bs{s} \in \{0,1\}^V}
\prod_{uv \in E} (1+ A_{uv} s_u s_{v} ) \prod_{v \in V} \exp (h_v s_v)  \\
& = \sum_{S \subset E}
\Big(
\prod_{uv' \in S} A_{uv'} \Big)
\prod_{v \in V} ( \sum_{s_v=0,1} s_v^{d_v(S)} \exp (h_v s_v))  \\
&= 
\prod_{v \in V} \exp (h_v )
\sum_{S \subset E, U \subset V \atop S \text{and} U \text{are not ``adjacent''}}
\prod_{uv \in S} A_{uv}
\prod_{v \in U} B_v \\
&=
\prod_{v \in V} \exp (h_v ) \quad 
p(G';\bs{A},\bs{B}),
\end{align}
where $d_v(S)$ is the number of edges in $S$ connecting to $v$,  
$A_{uv}={\rm e}^{J_{uv}}-1, B_v = {\rm e}^{- h_v}$ and $G'$ is a bipartite graph obtained by adding a new vertex
on each edge of $G$.
\end{proof}

\begin{remark}
{\rm
The Bethe approximation can also be applied to the computation of the
permanent of non-negative matrices.  From the combinatorial viewpoint,
this problem is related to the (weighted) perfect matching problem on
complete bipartite graphs.  Vontobel analyzed this problem and pose a
conjecture analogous to Eq.(\ref{eq:Zfinite}) \cite{vontobelbethe}. This
conjecture implies the inequality between the permanent and its Bethe
approximation, $Z \geq Z_B$, given his formula
Eq.(\ref{eq:PascalChar}).  The statement $Z \geq Z_B$ is, however,
directly proved by Gurvits, generalizing Schrijver’s permanental
inequality \cite{Gurvits}.
\rm}
\end{remark}

\bibliographystyle{unsrt}
\bibliography{BP}

\end{document}